\documentclass{amsart}
\usepackage{graphicx}
\newtheorem{theorem}{Theorem}
\newtheorem{lemma}[theorem]{Lemma}
\newtheorem{proposition}[theorem]{Proposition}
\newtheorem{corollary}[theorem]{Corollary}

\theoremstyle{remark}
\newtheorem{remark}[theorem]{Remark}

\def\lk{{\mathrm{lk}}}

\begin{document}

\title[SYMMETRY OF LINKS AND CLASSIFICATION OF LENS SPACES]
{SYMMETRY OF LINKS AND CLASSIFICATION OF \\LENS SPACES}

\author{J\'{o}zef H. Przytycki}
\address{Department of Mathematics, The George Washington University,
Washington, DC 20052, USA}
\email{przytyck@research.circ.gwu.edu}

\author{Akira Yasuhara}
\address{Department of Mathematics, Tokyo Gakugei University, 
Nukuikita 4-1-1, Koganei, Tokyo 184-8501, Japan \newline
{\em Current address}, {\rm October 1, 1999 to September 30, 2001}: 
Department of Mathematics, The George Washington University, 
Washington, DC 20052, USA}
\email{yasuhara@u-gakugei.ac.jp}


\begin{abstract}
We give a concise proof of a classification of lens spaces 
up to orientation-preserving homeomorphisms. 
The chief ingredient in our proof is a study of 
the Alexander polynomial of \lq symmetric' links in $S^3$. 
\end{abstract}

\maketitle
Let $T_1$ and $T_2$ be solid tori, and  
let $m_i$ and $l_i$ be the meridian and longitude of $T_i$ $(i=1,2)$. 
The {\em lens space} $L(p,q)$ is a 3-manifold that is 
obtained from $T_1$ and $T_2$ by identifying 
their boundaries in such a way that 
$m_2=p l_1+q m_1$ and 
$l_2=\overline{q} l_1+ r m_1$, where $(p,q)=1$ and $q\overline{q}-pr=1$.

In 1935, Reidemeister classified lens spaces 
up to orientation-preserving {\em PL} homeomorphisms \cite{Rei}. 
This classification was generalized to the topological category  
with the proof of the Hauptvermutung by 
Moise in 1952 \cite{Moi}. Meanwhile, Fox had outlined an aproach to  
classification up to homeomorphisms which would not require the 
Hauptvermutung; see \cite[Problem 2]{Eil}, \cite{Fox}. 
This was implemented later by Brody \cite{Bro}. 
We refer the reader to \cite{Gor} for history of 
classifications of lens spaces. 

In this paper, we give a concise proof of 
a classification of lens spaces up to orientation-preserving 
homeomorphisms. 
Our method is motivated by that of Fox-Brody. 
While the chief ingredient in their proof was a study of 
the Alexander polynomial of knots in lens spaces, 
we study the Alexander polynomial of \lq symmetric' links 
in $S^3$.

For an oriented 3-manifold $M$ with finite 
first homology group, the {\em linking form} 
\[\lk_M:H_1(M;{\Bbb Z})\times H_1(M;{\Bbb Z})
\longrightarrow{\Bbb Q}/{\Bbb Z}\]
is defined as follows \cite{Alex1}, \cite{Alex2}.  
Let $x$ and $y$ be $1$-cycles in $M$ that represent 
elements $[x]$ and $[y]$ of $H_1(M;{\Bbb Z})$ 
respectively. Suppose that 
$nx$ bounds a $2$-chain $c$ for some $n\in{\Bbb Z}$. Then 
\[\lk_M([x],[y])=\frac{c\cdot y}{n}\in{\Bbb Q}/{\Bbb Z},\]
where $c\cdot y$ is the intersection 
number of $c$ and $y$. 

Let $\Delta_K(t)$ be the 
{\em  Conway-normalized Alexander polynomial} of $K$, i.e., 
$\Delta_K(t)=\nabla_K(t^{-1/2}-t^{1/2})$, where 
$\nabla_K(z)$ is the Conway polynomial. 

\begin{theorem}
Let $\rho:S^3\longrightarrow L(p,q)$ be the $p$-fold cyclic cover
and $K$ a knot in $L(p,q)$ that represents a generator 
of $H_1(L(p,q);{\Bbb Z})$. If $\Delta_{\rho^{-1}(K)}(t)=1$, then 
$\lk_{L(p,q)}([K],[K])= q/p$ or $=\overline{q}/p$ in ${\Bbb Q}/{\Bbb Z}$. 
\end{theorem}

Before proving Theorem 1, we obtain the classification of lens spaces 
as its corollary. 

\begin{corollary} Two lens spaces 
$L(p,q)$ and $L(p,q')$ are equivalent up to
orientation-preserving homeomorphisms if and only if  
$q\equiv q'$ $($mod $p)$ or $qq'\equiv 1$ $($mod $p)$.  
\end{corollary}

\begin{proof}
The \lq if' part of the corollary is well-known and easy to see. 
It is enough to show the \lq only if' part. 
Since $pl_1+qm_1$ bounds a meridian disk $D_2$ of $T_2$, 
$pl_1$ bounds a 2-chain $qD_1\cup D_2$, 
where $qD_{1}$ is a disjoint union of $q$ copies of 
a meridian disk $D_1$ of $T_1$.  So we have  
$\lk_{L(p,q)}([l_1],[l_1])=
(qD_1\cup D_2)\cdot l_1/p=q/p$.  
Note that $\rho^{-1}(l_1)$ is the trivial knot, 
and hence $\Delta_{\rho^{-1}(l_1)}(t)=1$. 
Suppose that there is an orientation-preserving homeomorphism
$f:L(p,q)\longrightarrow L(p,q')$. Then 
$\lk_{L(p,q')}([f(l_1)],[f(l_1)])=\lk_{L(p,q)}([l_1],[l_1])$ and 
$\Delta_{\rho^{-1}(f(l_1))}(t)=\Delta_{\rho^{-1}(l_1)}(t)$. 
Therefore, by Theorem 1, we have $q/p\equiv q'/p$ or 
$\equiv\overline{q'}/p$ (mod $1$). This implies that  
$q\equiv q'$ (mod $p$) or $qq'\equiv 1$ (mod $p$).
\end{proof}

\begin{remark} 
Since $\lk_{L(p,q)}([nl_1],[nl_1])= n^2q/p$, 
the set $\{n^2q/p|1\leq n<p\}$ of rational numbers is an 
invariant of $L(p,q)$ up to orientation-preserving homotopy. 
A homotopy classification of lens spaces was obtained 
by Whitehead \cite{Whi}: $L(p,q)$ and $L(p,q')$ are equivalent up 
to orientation-preserving homotopy if and only if $qq'\equiv n^2$ 
(mod $p$) for some integer $n$.
The necessity of the condition is given by this invariant. 
\end{remark}

In order to prove Theorem 1, we need some preliminaries. 
In the following lemma, we consider the Alexander polynomial of
links with a certain kind of symmetry.

\begin{lemma} 
Let $r$ be a prime and $s$ a positive integer. 
Let $L$ be an $r^s$-periodic link in $S^3$ and 
let $L'$ be a link obtained from $L$ by changing a set of crossings that 
is the ${\Bbb Z}_{r^s}$-orbit of a single crossing in 
the periodic diagram of $L$. 
For an integer $q$, let $L(q)$ $($resp. $L'(q))$ be a link obtained 
from $L$ $($resp. $L')$ by 
adding $-q$-full twists as illustrated in Figure $1;$ equivalently we perform 
a $1/q$-surgery along the fixed point set of the 
${\Bbb Z}_{r^s}$-action. Then 
\[\Delta_{L(q)}(t)\equiv\Delta_{L'(q)}(t)\ {mod}\ 
(t^{r^s}-1,r).\]
\end{lemma}

\begin{proof}
In Figure 1(b), 
a tangle $T$ and the full-twists part are commutative.  
Therefore by arguments similar to that in the proof of 
Lemma 2.3 in \cite{Prz}, we obtain  
\[\Delta_{L(q)}(t)\equiv\Delta_{L'(q)}(t)\ {\mathrm {mod}}\ 
((t^{-1/2}-t^{1/2})^{r^s},r).\]
Since $(t^{-1/2}-t^{1/2})^{r^s}\equiv t^{-{r^s}/2}-t^{{r^s}/2}$ 
(mod $r$), we have the conclusion. 
\end{proof}

\begin{figure}[tb]
\includegraphics[trim=0mm 0mm 0mm 0mm, width=.75\linewidth]
{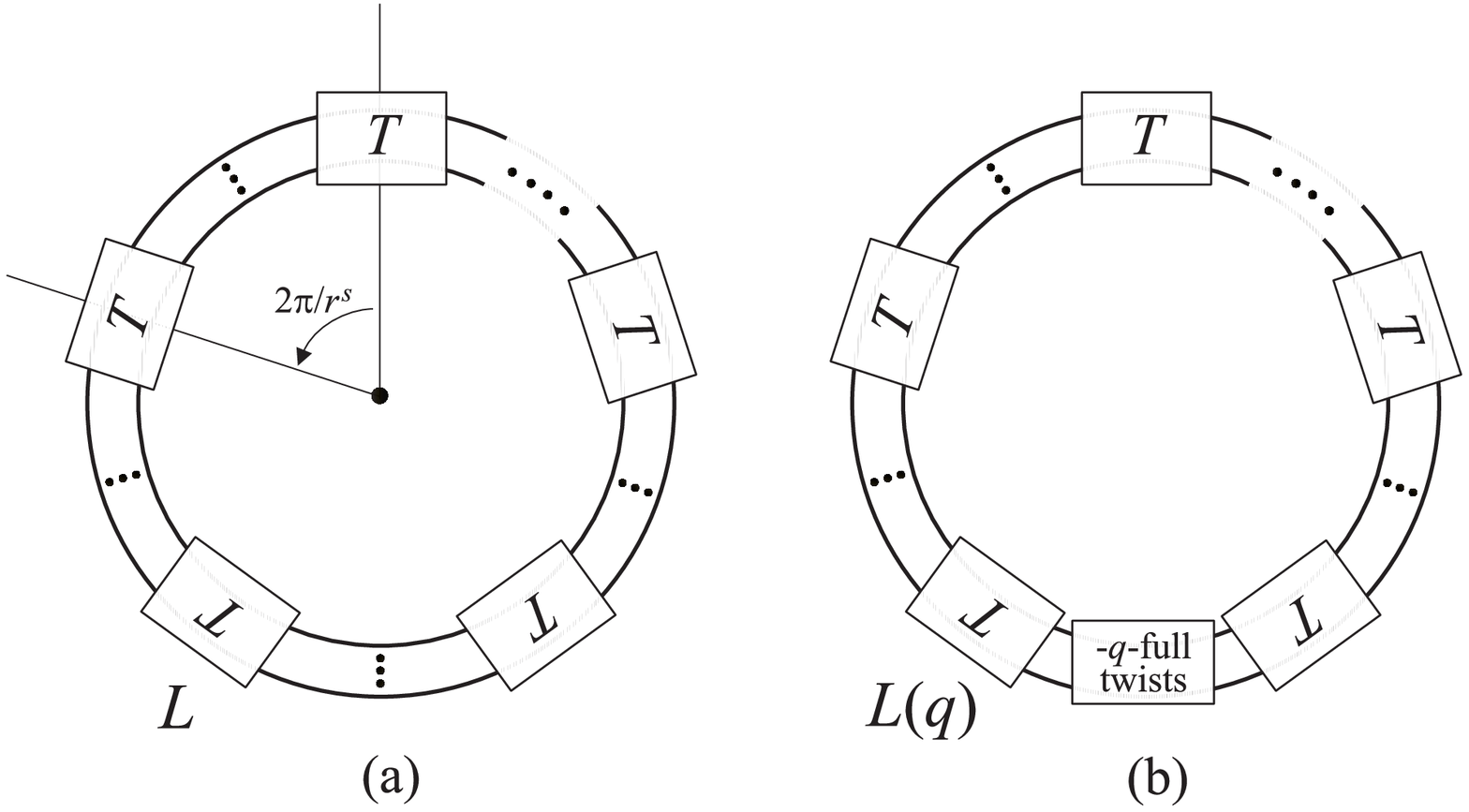}
\caption{}
\label{}
\end{figure}

A link in $S^3$, the universal cover of a lens space, that covers 
a knot in the lens space has the above-mentioned symmetry.  
By Lemma 4, we obtain the following proposition. 

\begin{proposition}
Let $\rho:S^3\longrightarrow L(p,q)$ be the $p$-fold cyclic cover 
and let $K$ and $K'$ be knots in $L(p,q)$. 
If $K$ and $K'$ are homologous in $L(p,q)$ 
and if $p$ is divisible by the $s$th power $r^s$ of a prime integer $r$, 
then 
\[\Delta_{\rho^{-1}(K)}(t)\equiv\Delta_{\rho^{-1}(K')}(t)\ {mod}\ 
(t^{r^s}-1,r).\]
\end{proposition}
  
\begin{proof} We may assume that both $K$ and $K'$ are contained in $T_1$. 
We note that the $p$-fold cover $S^3$ is obtained 
from $\tilde{T_1}$ and $\tilde{T_2}$ by identifying 
their boundaries in such a way that 
$\tilde{m_2}=\tilde{l_1}+q \tilde{m_1}$, where 
$\tilde{T_i}$ is a solid torus that is the $p$-fold cover 
of $T_i$, and $\tilde{m_i}$ and $\tilde{l_i}$ are the meridian and 
longitude of $\tilde{T_i}$ $(i=1,2)$. 
If two knots in $L(p,q)$ are homologous, then they are also homotopic. 
So they differ in a finite number of crossings. 
Thus it suffices to 
consider the case in which $K'$ is obtained from $K$ by changing a 
single crossing $c$. 
Then $\rho^{-1}(K)$ and $\rho^{-1}(K')$ are 
$p$-periodic in $\tilde{T_1}$, 
$\rho^{-1}(K')$ is obtained from $\rho^{-1}(K)$ by 
changing the crossings $\rho^{-1}(c)$, and 
a covering translation $\phi$ of $\tilde{T_1}$ 
generates the ${\Bbb Z}_p$-action on $\tilde{T_1}$. 
Set $p=ur^s$. 
Then we note that $\rho^{-1}(c)=\{c_{11},...,c_{1u}\}
\cup\cdots\cup\{c_{{r^s}1},...,c_{{r^s}u}\}$, where 
$\phi^{j-1}(c_{k1})=c_{kj}$ and  
$\phi^{u(k-1)}(c_{1j})=c_{kj}$ for $k=1,...,r^s,\ j=1,...,u$. 
This implies that there is a sequence 
$\rho^{-1}(K)=L_0,...,L_u=\rho^{-1}(K')$ of 
$r^s$-periodic links in $\tilde{T_1}$ such that 
$L_j$ is obtained from $L_{j-1}$ by changing the crossings 
$c_{1j},...,c_{{r^s}j}$ $(j=1,...,u)$. 
By Lemma 4, we have 
\[\Delta_{L_{j-1}}(t)\equiv\Delta_{L_j}(t)\ \mathrm{mod}\ 
(t^{r^s}-1,r).\]
This completes the proof. 
\end{proof}

\begin{proof}[Proof of Theorem $1$.] 
Let $K_n$ be a knot in $\partial T_1\subset L(p,q)$ that represents 
$n[l_1]+[m_1]\in H_1(L(p,q);{\Bbb Z})$. 
We may assume that $K$ is homologous to $K_n$ for some $n$ ($(n,p)=1$). 
Then $\rho^{-1}(K_n)$ is the $(n,p-qn)$-torus knot since 
$S^3$ is obtained from $\tilde{T_1}$ and $\tilde{T_2}$ by 
identifying their boundaries such that 
$\tilde{m_2}=\tilde{l_1}+q \tilde{m_1}$. 
It is well-known that 
\[\Delta_{\rho^{-1}(K_n)}=
t^{-g}\frac{(1-t)(1-t^{n(p-qn)})}{(1-t^n)(1-t^{p-qn})},\]
where $g=(n(p-qn)+1-(p-qn+n))/2$. 
By Proposition 5, 
\[t^{-g}
\frac{(1-t)(1-t^{n(p-qn)})}{(1-t^n)(1-t^{p-qn})}
\equiv 1 \ \mathrm{mod}\ 
(t^{r^s}-1,r),\]
where $r^s$ is a divisor of $p$ and $r$ is a prime integer. 
By elementary calculations, we obtain 
$n^2\equiv 1$ or $\equiv \overline{q}^2$ (mod $r^s$). 
So we have 
$n^2\equiv 1$ or $\equiv \overline{q}^2$ (mod $p$). 
Since $\lk_{L(p,q)}([K],[K])=n^2q/p$,  
we have the required result.  
\end{proof}

\bibliographystyle{amsplain}

\end{document}